\documentclass[12pt]{article}

\def\@typesizes{%
       \or{5}{6.5}\or{6}{7.5}\or{7}{8.5}\or{8}{11}\or{9}{12}%
       \or{10}{13}
       \or{\@xipt}{14}\or{\@xiipt}{15}\or{\@xivpt}{18}%
       \or{\@xviipt}{20}\or{\@xxpt}{24}}

\usepackage{amsthm}
\usepackage{amsmath}
\usepackage{amssymb}
\usepackage{graphicx}
\usepackage{enumerate}

\numberwithin{equation}{section}
\numberwithin{figure}{section}

\theoremstyle{plain}
\newtheorem{theorem}{ Theorem}[section]
\newtheorem{proposition}[theorem]{ Proposition}
\newtheorem{lemma}[theorem]{ Lemma}
\newtheorem{corollary}[theorem]{ Corollary}
\newtheorem{example}[theorem]{ Example}
\newtheorem{remark}[theorem]{ Remark}
\newtheorem{definition}[theorem]{ Definition}
\newtheorem{conjecture}{ Conjecture}

\def\BET{\begin{theorem}}
\def\ENT{\end{theorem}}
\def\BEP{\begin{proposition}}
\def\ENP{\end{proposition}}
\def\BEL{\begin{lemma}}
\def\ENL{\end{lemma}}
\def\BEC{\begin{corollary}}
\def\ENC{\end{corollary}}
\def\BEE{\begin{example} \rm}
\def\ENE{\end{example}}
\def\BER{\begin{remark} \rm}
\def\ENR{\end{remark}}
\def\BED{\begin{definition} \rm}
\def\END{\end{definition}}
\def\BECJ{\begin{conjecture}}
\def\ENCJ{\end{conjecture}}

\def\bea{\begin{eqnarray}}
\def\eea{\end{eqnarray}}
\def\beq{\begin{equation}}
\def\eeq{\end{equation}}
\def\beal{\begin{align*}}
\def\eeal{ \end{align*} }

\def\bfb{{\bf b}}
\def\bfc{{\bf c}}

\def\bfn{{\bf n}}

\def\bft{{\bf t}}

\def\bfB{{\bf B}}

\def\bfN{{\bf N}}

\def\bfT{{\bf T}}

\def\bbN{{\mathbb N}}

\def\bbR{{\mathbb R}}

\def\bbZ{{\mathbb Z}}

\def\cB{{\mathcal B}}

\def\cH{{\mathcal H}}

\newcommand{\ve}{\varepsilon}

\newcommand{\vt}{\vartheta}
\newcommand{\vph}{\varphi}

\begin{document}

\begin{center}
{\large{\bf Geometrically induced spectral effects in tubes
\\[.2em] with a mixed Dirichlet-Neumann boundary}}
\end{center}

\bigskip

\centerline{Fedor L. Bakharev$^{a,b}$ and Pavel Exner$^{c,d}$}

\bigskip

\begin{center}

\emph{a) St. Petersburg State University, Mathematics and Mechanics Faculty, 7/9 Universitetskaya nab., St. Petersburg, 199034 Russia} \\
\emph{b) St. Petersburg State University, Chebyshev Laboratory, 14th Line V.O., 29B, Saint Petersburg 199178 Russia} \\
\emph{c)Department of Theoretical Physics, Nuclear Physics Institute, Czech Academy of Sciences, 25068 \v Re\v z near Prague, Czechia} \\
\emph{d)Doppler Institute for Mathematical Physics and Applied Mathematics, Czech Technical University, B\v rehov\'a 7, 11519 Prague, Czechia} \\
\texttt{fbakharev@yandex.ru, f.bakharev@spbu.ru, exner@ujf.cas.cz}

\end{center}

\bigskip

{\small

\noindent {\bf Abstract}: We investigate spectral properties of the Laplacian
in $L^2(Q)$, where $Q$ is a tubular region in $\mathbb{R}^3$ of a
fixed cross section, and the boundary conditions combined a
Dirichlet and a Neumann part. We analyze two complementary
situations, when the tube is bent but not twisted, and secondly, it
is twisted but not bent. In the first case we derive sufficient
conditions for the presence and absence of the discrete spectrum
showing, roughly speaking, that they depend on the direction in
which the tube is bent. In the second case we show that a constant
twist raises the threshold of the essential spectrum and a local
slowndown of it gives rise to isolated eigenvalues. Furthermore, we
prove that the spectral threshold moves up also under a sufficiently
gentle periodic twist.

\medskip

\noindent {\bf Keywords}: Laplacian, Dirichlet-Neuman boundary, tube, discrete
spectrum

\medskip

\noindent {\bf MSC}: 81Q37, 35J05

}

\bigskip


\section{Introduction} 

Relations between spectral properties and geometry belong to
trademark topics in mathematical physics. A particularly interesting
class of problems concerns spectra of the Laplacians and related
operators in tubular regions which has various applications, among
others they are used to model waveguide effects in quantum systems.
The turning point here was the seminal observation that `bending
means binding', that is that the Dirichlet Laplacian in a tube of a
fixed cross section which is bent but asymptotically straight has a
nonempty discrete spectrum\footnote{Although it is not important in
the present three-dimensional context, we note this result extends
to tubes in higher dimensions \cite{ChDuFrKr}. In other situations
involving geometrically induced eigenvalues the effect may not be that
robust -- see, e.g. \cite{LoOu}.} -- see, e.g., \cite{DuEx}. It
inspired a long series of investigations, for a survey we refer to
the monograph \cite{ExKo15}
and the bibliography there.

A nontrivial geometry can be manifested not only in the shape of the
tube but also in the boundary conditions entering the definition of
the Laplacian. A simple but striking example can be found in
\cite{DiKr}: an infinite planar strip of constant width whose one
boundary is Dirichlet and the other Neumann exhibits a discrete
spectrum provided the Dirichlet boundary is bent `inward' while in
the opposite case the spectral threshold remains preserved. One is
naturally interested whether this effect has three-dimensional
analogue. The geometry is substantially richer in this case, of
course, nevertheless our first main result -- see Theorems~\ref{thm:exist}
and \ref{thm:nonexist} below -- provides an affirmative answer of a
sort to this question, namely that some bending
directions are favorable from the viewpoint of the discrete spectrum
existence and some are not.

Another class of geometric deformations are tube twistings. In
general, they act in the way opposite to bendings: to produce bound
states of the Dirichlet Laplacian supported by a locally twisted
tube of a non-circular cross section, an additional attractive
interaction must exceed some critical strength \cite{EkKoKr08}. On
the other hand, a discrete spectrum may arise in a tube which is
constantly twisted and the twist is locally slowed down
\cite{ExKo05}. Note that these results have a two-dimensional
analogue, namely a Hardy inequality in planar strips where the
Dirichlet and Neumann condition suddenly `switch sides'
\cite{KoKr08} and the appearance of a nontrivial discrete spectrum
when a sufficiently long purely Neumann segment is inserted in
between \cite{BorCar2,DiKr2}.

In the second part of the paper we examine twisted tubes with a
mixed Dirichlet-Neumann boundary. We show that the effect of
twisting and its local slowdowns is present again,
cf.~Proposition~\ref{Pr5.3}, now it may occur
also if the tube cross section itself exhibits a rotational symmetry
but the boundary conditions violate it. Furthermore, we consider a
wider class of tubes where the twist is not constant along the tube
but only periodic and ask whether in this case too the threshold of
the essential spectrum moves up; in Theorem~\ref{Th1} we demonstrate this property for twists
that are sufficiently gentle.

The main results of the paper indicated above, concerning the
effects of bending and twisting, constant and
periodic, are presented and proved in Sections~\ref{s:bending},
\ref{s:twisting}, and \ref{s:periodic},
respectively. Before coming to it, we collect in the next
two preliminary
sections the needed properties of the tubes and of the operators
involved.

\section{Preliminaries: geometry of the waveguide} \label{s:geometry}

Let us begin with a curve $\ell:\bbR \to \bbR^3$ that will play the
role of waveguide axis supposed to be a $C^3$-diffeomorphism of the
real axis $\bbR$ onto $\ell(\bbR)$. Without loss of generality we
may parametrize it by its arc length, that is, to assume that
$\dot{\ell_1}(z)^2+\dot{\ell_2}(z)^2+\dot{\ell_3}(z)^2=1$, where by
$\dot{\ell_j}$ we denote the derivative of function $\ell_j$ with
respect to the variable $z$. Dealing with the curve~$\ell$, we want
to associate with it the Frenet frame, i.e. the orthonormal triad of
smooth vector fields $\{\bft, \bfn, \bfb\}$ called respectively the
tangent, normal, and binormal vectors, defined as follows
$$ 
\bft=\dot{\ell}\,,\quad \bfn= \kappa^{-1} \ddot{\ell}\,, \quad
\bfb=\bft\times\bfn\,.
$$ 
where the cross denotes the vector product in $\bbR^3$ and
$\kappa:=|\ddot{\ell}|$ is the curvature of~$\ell$. Put like that, the
Frenet frame may not exist, in particular, because it is necessary
to assume that $\kappa>0$ holds to make sense of the definition of the
normal and binormal. If a part of $\ell$ is a straight line segment,
i.e. $\kappa=0$ holds on it identically, one can employ any fixed triad
one element of which coincides with the tangent vector. With a
slight abuse of terminology we will say that $\ell$ possesses a
\emph{global Frenet frame} if triads corresponding to its straight
and non-straight parts can glued together smoothly, modulo a
rotation of the Frenet parts on a fixed angle around the appropriate
tangent vector, see \cite{EkKoKr08, ExKo04} or
\cite[Sec.~1.3]{ExKo15}.

In such a case the Serret-Frenet formul{\ae} give
\begin{equation} \label{2.1-1} \dot{\bft}=\kappa \bfn\,,\quad
\dot{\bfn}=-\kappa\bft+\tau\bfb\,, \quad \dot{\bfb}=-\tau\bfn\,, \end{equation}
where $\tau$ stands for the torsion of the curve $\ell$. Given a
function $\beta\in C^1(\bbR)$ we introduce further a general moving
frame $\{\bfT_\beta, \bfN_\beta, \bfB_\beta\}$ by
\begin{equation} \label{2.2} \bfT_\beta=\bft\,,\quad \bfN_\beta= \bfn\cos \beta-
\bfb \sin \beta\,,\quad \bfB_\beta= \bfn\sin \beta+\bfb \cos
\beta\,; \end{equation}
the equations \eqref{2.1-1} show that this triad elements satisfy
the relations
\begin{equation} \label{2.1-2}
\begin{array}{l}
\dot{\bfT}_\beta=\kappa (\bfN_\beta \cos \beta +\bfB_\beta
\sin\beta)\,,\\[.3em]
\dot{\bfN}_\beta=-\kappa \bfT_\beta \cos\beta + (\tau-\dot{\beta})\bfB_\beta\,,
\\[.3em]
\dot{\bfB}_\beta=-\kappa \bfT_\beta \sin\beta
-(\tau-\dot{\beta})\bfN_\beta\,.
\end{array}
\end{equation}
A particular choice $\beta(z)=\int_{-\infty}^z \tau(s)ds$ yields the
so-called \emph{parallel transport frame}, in the physics literature
often referred to as the \emph{Tang frame}.

Let next $\omega$ be a two-dimensional bounded domain (by definition open
connected set) with the boundary $\partial\omega$ supposed to be piecewise $C^2$;
we suppose that $\omega$ has no cusps. The waveguide we are
going to consider is defined as the tube $Q_{\ell,\beta}$ obtained
by moving the cross-section $\omega$ along the reference curve
$\ell$ keeping its position fixed with respect to the frame
\eqref{2.2}. More precisely, we set
\begin{equation} \label{2.3} Q_{\ell,\beta} =\{ X\in\bbR^3 \colon
X=\ell(x_3)+x_1\bfN_\beta+x_2 \bfB_\beta, \, x'=(x_1,x_2)\in\omega,
\,z=x_3\in \bbR\}. \end{equation}
Denoting $a:=\sup_{x'\in \omega} |x'|$ and assuming that
\begin{equation} \label{a} a\|\kappa\|_\infty<1
\end{equation}
one can check easily that the formula \eqref{2.3} induces a local
$C^1$-diffeomorphism between $Q_{\ell,\beta}$ and the straight tube
$Q=\omega\times \bbR$. We will assume, without repeating it every
time, that this diffeomorphism is \emph{global}, in other words,
that the tube $Q_{\ell,\beta}$ has no self-intersections.

With the eye on the definition of the operators in the next section
we divide the boundary $\partial\omega$ into two parts. One denoted
as $\gamma_D$ is assumed to be a union of a finite number of arcs,
each of a positive measure, while its complement
$\partial\omega\setminus \gamma_D$ is denoted as $\gamma_N$. The
pair $(\omega,\gamma_D)$ is called \emph{rotationally invariant} if
$\omega$ and $\gamma_D$ are both rotationally invariant with respect
to the origin. From the viewpoint of this paper, this trivial case
that can occur only if $\omega$ is a disc, an annulus, centered at the origin,
and each connected component of
$\partial\omega$ is circle belonging to only one of the sets
$\gamma_D, \gamma_N$. In the following we will consider only
\emph{rotationally non-invariant} pairs $(\omega,\gamma_D)$.

Let us now specify two types of geometric deformations which we will
consider in this paper. We say that the tube $Q_{\ell, \beta}$ is
\emph{bent} if the reference curve~$\ell$ is not a straight line,
that is, the curvature $\kappa$ does not vanish identically.
Furthermore, the tube $Q_{\ell, \beta}$ is said to be \emph{twisted}
if the pair $(\omega, \gamma_D)$ is not rotationally invariant and
$\tau-\dot{\beta}\ne 0$. Looking at the equations \eqref{2.1-2} it
is obvious that these perturbations are mutually independent. For
the sake of simplicity, we will study separately the following two
cases:
\begin{itemize}
\item[(i)]
\emph{bending without twisting:} $\kappa\ne 0$ and
$\tau-\dot{\beta}=0\,$, and
\item[(ii)]
\emph{twisting without bending:} $\kappa = 0$ and $\tau-\dot{\beta}\ne
0\,$.
\end{itemize}
We will deal with the case (i) in Section~\ref{s:bending} and with (ii) in Sections~\ref{s:periodic} and \ref{s:twisting}, respectively.

\section{Preliminaries: definition of the operator} 

Before introducing the operator of our interest, we need a couple of
auxiliary notions. We denote by $\lambda_1$ the lowest eigenvalue of
the problem
\begin{equation} \label{psi} -\Delta' \psi =\lambda\psi \;\mbox{ in } \omega\,,
\quad \psi=0 \;\mbox{ on } \gamma_D\,, \quad \partial_n \psi=0
\;\mbox{ on } \gamma_N\,, \end{equation}
where $\Delta'=\nabla'\cdot\nabla'$ stands for the Laplace operator
with respect to the variables $x'$ and $\partial_n$ is the outward
normal derivative. The corresponding eigenfunction normalized in
$L^2(\omega)$ will be denote by $\psi_1$; we note that $\psi_1$ can
be chosen positive in $\omega$ and it satisfies the integral
identity
\begin{equation} \label{var} (\nabla' \psi_1,\nabla' \phi)_{\omega}=\lambda_1
(\psi_1, \phi)_{\omega} \qquad \forall \phi\in H^1_0(\omega, \gamma_D)\,,
\end{equation}
where $(\cdot,\cdot)_{\omega}$ is the natural scalar product in the
Lebesgue space $L^2(\omega)$ and $H^1_0(\omega,\gamma_D)$ consists
of functions from the first Sobolev space that vanish on $\gamma_D$.

Note that the continuity of $\psi_1$ on $\overline{\omega}$ is
non-trivial at the points of changing of boundary conditions. Such
a problem is studied, e.g., in the book \cite{NaPl}.

\BER
The crucial role in Section \ref{s:periodic} will play the following
condition
\begin{equation} \label{assump}
\partial_\vph \psi_1\ne 0 \quad \Leftrightarrow \quad
\int_\omega |\partial_\vph \psi_1|^2 \,\mathrm{d}x'\ne 0\,.
\end{equation}
If one considers Dirichlet boundary conditions only the validity of \eqref{assump}
can be achieved by requiring a rotational non-invariance of $(\omega, \gamma_D)$
-- see, e.g., \cite{Kr} or \cite{ExKo15} -- but if we admit mixed boundary
conditions it is not the case. For instance, let us consider Bessel function
$J_0$ of first kind and let $\nu_1$ and $\nu_2$ be the first two roots of
$J_0$. Then $J_0$ becomes the first eigenfunction of the Laplacian in any angular
sector of the annulus $\mathbb{B}_{\nu_1,\nu_2}=\{x'\colon |x'|\in(\nu_1,\nu_2)\}$
with Dirichlet conditions on the circular part of the boundary and Neumann condition on the radial
part, but at the same time we have $\partial_\varphi J_0=0$.
\end{remark}

The main object of our interest is the Laplace operator
$\widetilde{T}_{\ell,\beta}$ on $L^2(Q_{\ell,\beta})$ with mixed
Dirichlet-Neumann boundary conditions that can be associated with
the closed quadratic form
$$
\widetilde{a}_{\ell,\beta}[u]=\int_{Q_{\ell,\beta}}|\nabla_X
u|^2\,\mathrm{d}X, \quad u\in
H^1_0(Q_{\ell,\beta},\Gamma_{\ell,\beta})\,,
$$
where $\Gamma_{\ell,\beta}:=\{X\in\bbR^3 \colon
X=\ell(x_3)+x_1\bfN_\beta+x_2 \bfB_\beta, \, x'\in\gamma_D, \,
x_3\in\bbR\}$.

The diffeomorphism \eqref{2.3} can be used to map
$\widetilde{T}_{\ell,\beta}$ in the usual way \cite[Secs.~1.3 and
1.7]{ExKo15} to an operator on the straight tube $Q$ in which the
geometry is encoded in the coefficients. Specifically, in case (i)
the operator $\widetilde{T}_{\ell,\beta}$ is unitarily equivalent to
operator $T_{\ell, \beta}$ associated with the quadratic form
\begin{equation} \label{form(i)}
a_{\ell,\beta}[u]:=\int_Q \big[g (|\partial_1 u|^2+|\partial_2
u|^2)+g^{-1} |\partial_3 u|^2\big] \,\mathrm{d}x\,, \quad u\in
H^1_0(Q,\Gamma,g)\,,
\end{equation}
in the weighted Lebesgue space $L^2(Q,g)$ with the scalar product
$$
(u,v)_{Q,g}=(u g, v)_{Q}\,,
$$
where $g(x):=1-(x_1 \cos \beta(x_3)+x_2 \sin\beta(x_3))\kappa(x_3)$; in \eqref{form(i)} and in the following
we use $\Gamma$ as a shorthand for the set $\gamma_D\times \bbR$. In
particular, when the bending is absent, $\kappa=0$, and the
cross-section remain fixed in the parallel transport frame,
$\dot{\beta}=\tau$, the spectrum of $T_{\ell,\beta}$ is found by
separation of variables: it is purely continuous and equal to
$[\lambda_1,+\infty)$, where $\lambda_1$ is the eigenvalue appearing
in \eqref{var}. The question we address in the next section is under
which circumstances can a bending of such a tube give rise to a
nonempty discrete spectrum below the threshold $\lambda_1$.

Similarly, in case (ii) the operator $\widetilde{T}_{\ell,\beta}$ is
unitarily equivalent to the operator $T_\beta$ associated with the
quadratic form
\begin{equation} \label{form(ii)}
a_\beta[u]:=\int_Q \big[|\partial_1 u|^2+|\partial_2
u|^2+|(\partial_3+\dot{\beta} \partial_\vph)
u|^2\big]\,\mathrm{d}x\,, \quad u\in H^1_0(Q,\Gamma_D)\,,
\end{equation}
in the space $L^2(Q)$. Here the role of an unperturbed system will
be played by tubes with a constant twisting; in
Sec.~\ref{s:twisting} below we will discuss what happens if the
twisting rate is modified locally.

\section{The effect of bending} \label{s:bending}

Let us first focus on spectral properties of $T_{\ell, \beta}$ if
the tube exhibits a bending without twisting, i.e. the case (i)
indicated in Sec.~\ref{s:geometry}. To state the results, we need to
introduce two quantities,
$$
A_j= \frac12 \int_{\partial\omega} n_j |\psi_1|^2 \,\mathrm{d S}(x')=
\int_{\omega} \psi_1 \partial_j \psi_1 \,\mathrm{d}x'\,,\quad j=1,2\,,
$$
where $(n_1,n_2)$ are the components of outward unit normal to the
boundary $\partial \omega$. Note that while we use modulus in the
first expressions as it is common in quantum mechanics, we suppose
that the function $\psi_1$ is real-valued. This may be done without
loss of generality, since the operator commutes with the complex
conjugation, or in physical terms, the system is invariant with
respect to the time reversal.

We emphasise that the technique used in this section is not new.
The proof of Theorem~\ref{thm:exist} repeats the proof from \cite{DiKr}
almost literally. On the other hand, we want to mention that the proof of
Theorem~\ref{thm:nonexist} given below is much simpler than the proof of the analogous
result in \cite{DiKr}.


\begin{theorem}
\label{thm:exist}
If there exists a compact interval $I\subset \bbR$ such that
$$
\int_{I} \kappa(x_3) (A_1 \cos \beta(x_3)+A_2 \sin\beta(x_3))
\,\mathrm{d}x_3 <0
$$
and
\begin{equation}
\label{pointwise}
\kappa(x_3) (A_1 \cos \beta(x_3)+A_2 \sin\beta(x_3))\leq 0
\end{equation}
holds
for all $x_3\in \bbR\setminus I$, then
$$
\inf \sigma(T_{\ell,\beta})<\lambda_1\,.
$$
In particular, if the curvature $\kappa$ is compactly supported, the
discrete spectrum of $T_{\ell, \beta}$ is nonempty.
\end{theorem}

\noindent{\bf Proof.}
It is sufficient to find a trial function $u\in
H^1_0(Q,\Gamma,g)$ such that
\begin{equation} \label{trial} a_{\ell,\beta}[u]- \lambda_1\|u\|_{Q,g}^2<0\,,
\end{equation}
where $a_{\ell,\beta}$ is the quadratic form \eqref{form(i)}. We
will seek it in the form $u(x)=v(x_3) \psi_1(x')$, where $v$ is a
smooth function with compact support such that $v(x_3)=1$ for
$x_3\in I$ and
\begin{equation}
\label{delta} \|\partial_3 v\|_{L^2(\bbR)}=\delta\,,
\end{equation}
where $\delta$ is a parameter to be chosen. The assumption \eqref{a}
in combination with relation \eqref{delta} imply that there is a $C>0$
such that
$$
\int_Q g^{-1}|\partial_3 u|^2 \,\mathrm{d}x \leq C \delta^2\,,
$$
At the same time, the remaining part of the quadratic form in
question is
\begin{multline*}
\int_Q g (x) (|\nabla' u(x)|^2 -\lambda_1|u(x)|^2) \,\mathrm{d}x =\\
=-\int_\bbR  \kappa(x_3)\cos \beta(x_3) |v(x_3)|^2 \,\mathrm{d}x_3
\int_\omega x_1 (|\nabla' \psi_1(x')|^2 -\lambda_1|\psi_1(x')|^2)\,\mathrm{d}x'-\\
- \int_\bbR  \kappa(x_3)\sin \beta(x_3) |v(x_3)|^2 \,\mathrm{d}x_3
\int_\omega x_2 (|\nabla' \psi_1(x')|^2
-\lambda_1|\psi_1(x')|^2)\,\mathrm{d}x'\,,
\end{multline*}
where we have employed the explicit formula of $g(x)$ and relation
\eqref{var}. Next we note that
$$
A_j=-\int_\omega x_j (|\nabla' \psi_1(x')|^2
-\lambda_1|\psi_1(x')|^2)\,\mathrm{d}x'\,, \quad j=1,2\,.
$$
Indeed, let us insert $\phi(x')=x_j \psi_1(x')$ into the integral
identity \eqref{var} and rewrite it in the following way
\begin{align*}
-\int_\omega x_j (|\nabla' \psi_1(x')|^2
-\lambda_1|\psi_1(x')|^2)\,\mathrm{d}x'
&= \int_ {\omega} \psi_1(x') \nabla' x_j \cdot \nabla'\psi_1(x')\,\mathrm{d}x'\\
&= \int_{\omega} \psi_1(x')\partial_j\psi_1(x')\,\mathrm{d}x'\,.
\end{align*}
Under the stated assumptions the whole expression is negative and
choosing $\delta$ small enough we can make relation \eqref{trial}
satisfied. We know that for the tube without bending, $\kappa=0$, the
spectrum is purely essential and equal to $[\lambda_1,+\infty)$ and
it is easy to check that it remains preserved under a compactly
supported perturbation, hence the spectrum below $\lambda_1$ must be
in such a case discrete and nonempty.
\hfill$\Box$

\BER
Let us note that we defined the curvature $\kappa$ as a non-negative
function, an therefore the pointwise estimate~\eqref{pointwise} means that
$$
A_1 \cos \beta(x_3)+A_2 \sin\beta(x_3)\leq 0 \quad\mbox{if}\quad
\kappa(x_3)\ne 0\,.
$$
\end{remark}

On the other hand, we can state a condition under which the bending
does not move $\inf\sigma(T_{\ell,\beta})$ down, which means, in
particular, the absence of eigenvalues for (non-twisted) tubes with
a compactly supported curvature.

\begin{theorem} \label{thm:nonexist}
If for all $x_3\in \bbR$ and $x'\in \omega $ the inequality
\begin{equation} \label{pr23}
\kappa(x_3) (\partial_1 \psi_1(x') \cos \beta(x_3)+
\partial_2 \psi_1(x')\sin \beta(x_3)) \geq 0
\end{equation}
is valid, we have $\inf \sigma(T_{\ell,\beta})\geq \lambda_1$.
\end{theorem}
\noindent{\bf Proof.}
Fix an arbitrary function $u\in C_0^\infty(Q\cup (\partial Q\setminus \Gamma))$.
The function $\psi_1$ is strictly positive in $\omega$ as well as on $\gamma_N$.
The last assertion follows from Zaremba--Hopf--Oleinik lemma -- see, e.g., \cite{NaAI}
and references therein.
Then we can write it as $u(x)=\psi_1(x')v(x)$ with
some $v\in H^1(Q)$. The `shifted' quadratic form entering
\eqref{trial} can be estimated from below by neglecting the
non-negative term containing the derivative with respect to the
longitudinal variable $x_3$,
\begin{equation} \label{prop23}
a_{\ell,\beta}[u]-\lambda_1\|u\|_{Q, g}^2\geq
\int_Q g(x)\left(|\nabla' u(x)|^2-\lambda_1|u(x)|^2\right)\,\mathrm{d}x.
\end{equation}
The above described factorization yields the formula
$$
|\nabla'(\psi_1 v)|^2= |\psi_1|^2 |\nabla' v|^2 +\nabla' \psi_1
\cdot \nabla' (\psi_1 v^2)\,,
$$
which allows us to split the last integral in \eqref{prop23} into
two parts,
$$
J_1=\int_Q g(x)|\psi_1(x')|^2 |\nabla' v(x)|^2 \mathrm{d}x \geq 0
$$
and
$$
J_2=\int_Q g(x) \big(\nabla'\psi_1(x') \cdot\nabla' (\psi_1(x')
v^2(x))-\lambda_ 1|\psi_1(x')|^2 |v(x)|^2\big)\,\mathrm{d}x\,.
$$
Using the explicit form of $g(x)$ in combination with \eqref{var} we
get
\begin{multline*}
\hspace{-1em} J_2=-\int_Q x_1\cos\beta(x_3) \kappa(x_3) \big(\nabla'
\psi_1(x') \cdot\nabla' (\psi_1( x') v^2(x))-\lambda_1|\psi_1(x')|^2
|v(x)|^2\big)\,\mathrm{d}x \\
-\int_Q x_2\sin\beta(x_3) \kappa(x_3) \big(\nabla' \psi_1(x') \cdot\nabla'
(\psi_1(x') v^2(x))-\lambda_1|\psi_1(x')|^2 |v(x)|^2\big)\,\mathrm{d}x.
\end{multline*}
Inserting next the function $x_j\psi_1(x') |v(x',x_3)|^2$ into the
integral identity \eqref{var} as a test function with $x_ 3$ as
parameter, we obtain
\begin{multline*}
\int_{\omega} x_j \nabla'\psi_1(x')\cdot \nabla'(\psi_1(x')
v^2(x))\mathrm{d}x'-\lambda_1
\int_{\omega} x_j|\psi_1(x')|^2 |v(x)|^2\mathrm{d}x'=\\
=-\int_{\omega}\psi_1(x')v^2(x)  \nabla' x_j\cdot \nabla'
\psi_1(x')\, \mathrm{d}x' =-\int_{\omega}\psi_1(x') v^2(x) \partial_j
\psi_1(x') \,\mathrm{d}x'\,,
\end{multline*}
and consequently,
$$
J_2=\int_Q \kappa(x_3) (\partial_1 \psi_1(x') \cos \beta(x_3)+
\partial_2 \psi_1(x' )\sin \beta(x_3)) \psi_1(x') v^2(x)
\,\mathrm{d}x\,,
$$
which together the assumption \eqref{pr23} shows that
$a_{\ell,\beta}[u]- \lambda_1\|u\|_{Q,g}^2 \ge 0$ holds for any
$u\in C_0^\infty(Q\cup (\partial Q\setminus \Gamma))$. To finish the proof it is enough to observe
that this set is dense in $H^1_0(Q,\Gamma)$.
\hfill$\Box$

\medskip

In the particular case where the waveguide axis is a planar curve
and the Jacobian depends only on the coordinates $x_1$ and $x_3$, in
other words, $\tau(x_3)=\beta(x_3)=0$ holds for all $x_3\in \bbR$,
we have $\bfN=\bfn$, $\bfB=\bfb$ and the quadratic form expression
simplifies to
$$
a_{\ell,0}[u]=\int_Q (1-x_1 \kappa(x_3))^{-1} |\partial_3 u(x)|^2
+(1-x_1 \kappa(x_3))(|\nabla' u(x)|^2)\,\mathrm{d}x\,;
$$
then the above theorems lead to the following conclusions:
\begin{corollary}
If $A_1 <0$ and $\kappa\ne 0$ holds on a set of positive measure, then
$\inf \sigma(T_{\ell,0})<\lambda_1$.
\end{corollary}
\begin{corollary}
If for all $x'\in \omega$ the inequality
\begin{equation}
\label{restrictive}
\partial_1 \psi_1(x') \geq 0
\end{equation}
is valid, then $\inf \sigma (T_{\ell,0}) \geq \lambda_1$.
\end{corollary}
\BER
{\rm We note that when the curve is planar it is more natural
to consider the signed curvature. The corresponding slightly
modified reasoning is useful when comparing the conclusion with
the results of \cite{DiKr} as we are going to do below.}
\end{remark}
\BER
{\rm The condition \eqref{restrictive} seems to be rather restrictive,
however, one cannot weaken it as can be seen using an example constructed
as in the beginning of Sec.~3 in \cite{DiKr}, with the curvature support
consisting of two disjoint intervals. On the other hand, one may wonder
whether the assumption \eqref{restrictive} is not empty. This is not
the case, indeed, the simplest example is the square
$\omega=(0,1)\times(0,1)$ with $\gamma_D=\{x'\colon x_1=0\}$,
another example is the disc $\omega=\{x'\colon |x'|<1\}$ with
$\gamma_D=\{x'\colon |x'|=1, x_1<0\}$.}
\end{remark}

With the motivation explained in the introduction in mind, it
is instructive to compare the above results with the spectral
properties of a bent and asymptotically straight two-dimensional
guide which has one boundary Dirichlet and one Neumann. As we
have mentioned, it is known \cite{DiKr} that the discrete
spectrum of such a system is nonvoid if the total bending  of the
strip is positive (in fact, nonnegative) and the Dirichlet boundary
faces `inward', and on the other hand, the spectral threshold
remains preserved if the bend is simple, i.e. the curvature does not
change sign, and the Dirichlet condition is on the `outward' side.
This nicely fits with the above corollaries if we realize that the
normal, in other words, the $x_1$ direction points `'inwards' in a
bend and the outward normal derivative of $\psi_1$ is zero at the
Neumann segment(s) of the boundary and negative at the Dirichlet
one(s). Since cross sections with a piecewise smooth boundary are
covered by our assumptions, the simplest example is represented
by a wave\-guide with a rectangular cross-section and two flat sides.
If the bent sides of such a rectangular tube are
Dirichlet and Neumann and a fixed condition is chosen on each of the
flat sides, by separation of variables we get a direct
correspondence between the said two-dimensional properties and the
results obtained here.

\section{Twisting without bending} \label{s:twisting}

Let us now turn to the second class of geometric perturbations
indicated in Sec.~\ref{s:geometry} and discuss the situation when
the waveguide with a mixed Dirichlet-Neumann boundary is twisted
waveguide. Recall first how the situation looks like for tubes with
purely Dirichlet boundary. If the twisting is only local there it
does not affect the essential spectrum of the Laplacian, and
moreover, it does stabilize it against negative perturbations --
see, e.g., \cite{BoMaTr07, EkKoKr08}. This has consequences such as
the absence of weakly coupled bound states of Schr\"odinger
operators in twisted waveguides \cite{EkKoKr08, Kr}. If the twist is
not local but constant, then it even increases the threshold of the
essential spectrum of the Laplacian, and moreover, any local
slowdown of the constant twisting rate induces at least one bound
state of the corresponding operator \cite{ExKo05}.

Our aim in this section is to show that the behavior of twisted
tubes with a mixed Dirichlet-Neumann boundary is similar. We thus
suppose that $\dot{\beta}_0(x_3)=\alpha$ and introduce a model
eigenvalue problem on the cross section,
\begin{eqnarray}
&&-\Delta \psi^\alpha-\alpha^2 \partial_\vph^2\psi^\alpha=\lambda
\psi^\alpha \mbox{ in }\omega\,, \quad \psi^\alpha=0 \mbox{ on }
\gamma_D\,,\nonumber\\
\label{conormal}
&&\partial_n \psi^\alpha+\alpha^2(n_2x_1-n_1x_2)\partial_\varphi \psi^\alpha=0
\mbox{ on } \gamma_N\,,
\end{eqnarray}
where the expression in the left-hand side of \eqref{conormal} is the
co-normal derivative of the operator $(-\Delta-\alpha^2\partial^2_{\varphi})$.
We denote the smallest eigenvalue of this problem by
$\lambda_1^\alpha$, the corresponding normalized eigenfunction in
$L^2(\omega)$ will be $\psi_1^\alpha$; we note that $\psi_1^\alpha$
can be supposed to be positive without loss of generality and that
it satisfies the integral identity
$$
(\nabla' \psi_1^\alpha,\nabla' \phi)_{\omega}+\alpha^2 (\partial_\vph
\psi_1^\alpha,\partial_\vph \phi)_{\omega} =\lambda_1^\alpha
(\psi_1^\alpha, \phi)_{\omega} \qquad \forall \phi\in H^1_0(\omega,
\gamma_D)\,.
$$
\begin{proposition} 
If $\dot{\beta}_0=\alpha$ is constant the spectrum of
the positive self-adjoint operator $T_{\beta_0}$ coincides with the
interval $ [\lambda_1^\alpha,+\infty)$.
\end{proposition}
\noindent{\bf Proof.}
Similarly to the proof of Theorem~\ref{thm:nonexist}, we consider
function from a core of $T_{\beta_0}$ writing them as
$u(x)=\psi_1^\alpha(x') v(x)$ with $v\in H^1(Q)$, then a direct
calculation shows that
$$
a_{\beta_0}[u] - \lambda_1^\alpha \|u\|_{Q}^2 = \int_Q
|\psi_1^\alpha(x')|^2 \big(|\nabla' v(x)|^2+|(\partial_3 +\alpha
\partial_\vph)v(x)|^2\big)\,\mathrm{d}x\geq 0
$$
and by the density argument we obtain the estimate $\inf
\sigma(T_{\beta_0})\geq \lambda_1^\alpha$. To complete the proof it
is sufficient to construct in the standard way Weyl sequences for
any $\lambda \in [\lambda_1^\alpha,+\infty)$.
\hfill$\Box$

In this way a constant twisting, $\dot{\beta}_0=\alpha$, changes the
essential spectrum in a way depending on $\alpha$. Consider next a
local slowdown of the twist. Let $\theta=\theta(x_3)$ be a
$C^1$-function supported in a compact interval $I$ and assume that
the rotation angle $\beta$ is of the form
\begin{equation} \label{slow} \dot{\beta} (x_3)=\alpha-\theta(x_3)\,. \end{equation}
From the compactness of $\mathop{\rm supp} \theta$ it follows by a standard
perturbation argument that
$$
\inf \sigma_\mathrm{ess}(T_\beta)=\inf
\sigma(T_{\beta_0})=\lambda_1^\alpha\,.
$$

\medskip

\begin{proposition}
\label{Pr5.3}
Let $\dot{\beta}$ be given by the formula \eqref{slow} and
$$
\int_\bbR (|\dot{\beta}(x_3)|^2-\alpha^2)\,\mathrm{d}x_3 <0\,,
$$
then $\inf \sigma(T_\beta)<\lambda_1^\alpha$, and consequently,
$\sigma_\mathrm{disc}(T_\beta)\ne \varnothing$.
\end{proposition}
\noindent{\bf Proof.}
We employ trial functions $u\in H^1_0(Q, \Gamma_D)$ of the
factorized form $u(x)=\psi_1^\alpha(x') v(x_3)$, where $v$ is a
smooth function such that $v(x_3)=1$ if $x_3\in I$; it is easy to
check that for any $\delta>0$ one can choose $v$ so that $\|\partial_3
v\|_{L^2(\bbR)} =\delta$. A straightforward calculation then yields
$$
a_{\beta}[u]-\lambda_1^\alpha\|u\|_Q^2=\delta^2 +\int_\bbR
(|\dot\beta(x_3)|^2-\alpha^2)\,\mathrm{d}x_3 \int_\omega
|\partial_\vph\psi_1^\alpha(x')|^2 \,\mathrm{d}x'.
$$
Taking into account relation \eqref{assump} we find that for
$\delta$ small enough we have $a_{\beta}[u]
-\lambda_1^\alpha\|u\|_{Q}^2<0$ which completes the proof.
\hfill$\Box$

\section{Periodic twisting} \label{s:periodic}

Now we are going to consider a more general situation, bending still
absent and the twisting is non-constant but periodic leading to a
band-gap structure of the spectrum. It is natural to expect that a
higher twisting rate could increase the spectral threshold. Our aim
here is to demonstrate that it is indeed the case provided the
twisting is gentle. To state the result let us denote by $\beta\in
C^2(\bbR)$ the twisting function with a positive and 1-periodic
derivative $\dot{\beta}$. Let $\lambda_\dagger(\beta)$ be the
spectral threshold of the Laplacian with the mixed Dirichlet-Neumann
boundary conditions in the twisted tube $Q_\beta$. The main result
of the present section is the following theorem:
\begin{theorem} \label{Th1} Let pair $(\omega,\gamma_D)$ satisfies the assumption
\eqref{assump} and let $\vt_1, \vt_2\in C^2(\bbR)$ with $\dot{\vt}_1$,
$\dot{\vt}_2$ being positive 1-periodic functions and
$$
\int_0^1|\dot{\vt}_1|^2\, \mathrm{d}x_3 <\int_0^1|\dot{\vt}_2|^2\,
\mathrm{d}x_3\,,
$$
then there exists an $\ve_0>0$ such that the inequality
$\lambda_\dagger(\ve{\vt}_1) <\lambda_\dagger(\ve{\vt}_2)$ holds for
all $\ve\in(0,\ve_0)$. \end{theorem}

\noindent We will prove the theorem in several steps.

\subsection{Formulation of the problem}

If the function $\dot{\beta}$ is 1-periodic the spectrum of positive
self-adjoint operator $T_\beta$ is known to be purely essential
having a band-gap structure,
\begin{equation} \label{5.4} \sigma(T_\beta)=\bigcup_n\, \cB_n(\beta)\,. \end{equation}
One naturally expects it to be absolutely continuous, however, this
property has been so far established in some cases only
\cite[Chap.~9]{ExKo15}. Among spectral properties of periodic
waveguides, the gap opening for Laplace operator with various
boundary conditions has been discussed in many papers -- cf., e.g.,
\cite{BaNaRu, CaNaPe, Na, Yo98} -- in particular, for the case of
Dirichlet Laplacian in periodically twisted waveguide see
\cite{BoPa}.

To study the spectrum \eqref{5.4} we use Floquet--Bloch theory
(see, e.g., \cite{RS4}) and
decompose the operator $T_\beta$ into a direct integral of the
operator family $\{T_\beta(\eta)\}$ parametrized by the
quasi-momentum $\eta\in[-\pi,\pi]$. The fiber operators
$T_\beta(\eta)$ can be defined through their quadratic forms
$$
a_{\beta,\eta}[U]=\int_\Omega |\nabla'
U|^2+|(\partial_3+\dot{\beta}\partial_\vph)U|^2\,\,\mathrm{d}x\,,
\quad U\in \cH_\eta\,,
$$
related as usual to the corresponding sesquilinear forms
$A_{\beta,\eta}$ by $a_{\beta,\eta}[U] = A_{\beta,\eta}(U,U)$, where
$\Omega=\omega\times (0,1)$ is the periodicity cell and the form
domain $\cH_\eta$ consists of functions $U\in
H^1_0(\Omega,\gamma_D\times (0,1))$ satisfying a quasi-periodicity
condition
$$
U(x',1)=e^{i\eta}U(x',0)\,, \quad x'\in\omega\,.
$$
The inner product in $\cH_\eta$ is given by $\langle U,V\rangle =
(\nabla U, \nabla V)_\Omega$. It is not difficult to check that the
quadratic form $a_{\beta,\eta}$ is positive and closed, and
therefore associated with a unique self-adjoint operator
$T_\beta(\eta)$. Due to the compactness of the periodicity cell the
latter has a compact resolvent, and as a consequence, the spectrum
of operator $T_\beta(\eta)$ is purely discrete, in other words, a
sequence of eigenvalues
\begin{equation} \label{8}
0<\Lambda_{1,\beta}(\eta)\leq\Lambda_{2,\beta}(\eta)\leq\Lambda_{3,\beta}(\eta)\leq\ldots
\end{equation}
accumulating only at infinity; as the inequalities \eqref{8}
suggest, except the first some on them may not be simple. We denote
the corresponding eigenfunctions by $U_{k,\beta}\in \cH_\eta$, where
for the sake of simplicity the dependence on $\eta$ will usually not
be shown. We can choose them so that they satisfy the orthogonality
property
$$
(U_{j,\beta},U_{k,\beta})_\Omega=\delta_{j,k}\,,\quad
j,k=1,2,3,\ldots\,,
$$
where $\delta_{j,k}$ is the Kronecker symbol. The band functions
$\eta\mapsto \Lambda_{j,\beta}$ are known to be continuous and
$2\pi$-periodic so the sets
$$
\cB_j(\beta)=\{\Lambda_{j,\beta}(\eta)|\: \eta\in
[-\pi,\pi]\}\subset [0,+\infty)\,,
$$
the spectra bands, are closed intervals. In this notation the
overall spectral threshold can be written as
$$
\lambda_\dagger(\beta)=\inf \sigma(T_\beta)=\inf_{\eta\in[-\pi,\pi]} \Lambda_{1,\beta}(\eta)\,.
$$
In the absence of twisting, $\beta=0$, the eigenvalues and
eigenfunctions of the fiber operator are easily found explicitly,
\begin{equation} \label{10} V_{j,k}(x,\eta)=e^{i\eta x_3} e^{2\pi i k
x_3}\psi_j(x')\,, \,\, M_{j,k}(\eta)=(\eta+2\pi k)^2+\lambda_j\,,\,\, k\in\bbZ,\, j\in\bbN,
\end{equation}
where $(\psi_j,\lambda_j)$ is the $j$th eigenpair of the problem
\eqref{psi}. The family $\{\psi_j\}_{j=1}^\infty$ can be chosen to
be othonormal in $L^2(\omega)$ and the eigenvalue sequence
$\{\lambda_j\}_{j=1}^\infty$ is conventionally ordered in the
non-decreasing way counting multiplicities. Rearranging the sequence
$\{M_{j,k}(\eta)\}$ in the ascending order we obtain \eqref{8} for
$\beta=0$, in particular, $\Lambda_{1,0}(\eta)=\lambda_1+\eta^2$; we
note that these eigenvalues are simple unless $\eta=\pm \pi$.

\subsection{Small twisting, simple estimates}

Let us introduce a positive parameter $\ve\in (0,1)$ and discuss the
properties of the spectrum $\sigma(T_{\ve \beta})$ as $\ve\to 0$ for
a given function $\beta$. We start with the following simple lemma.
\begin{lemma} \label{Lemma 6.1} There is a constant  $C_\beta$ such that for
all $\ve\in (0,1)$, $\,\eta\in [-\pi,\pi]$, and $U\in \cH_\eta$ the
estimate
\begin{equation} \label{L1}
|a_{\ve\beta,\eta}[U]-a_{0,\eta}[U]|\leq C_\beta\, \ve \|\nabla
U\|_\Omega^2
\end{equation}
holds. \end{lemma}
\noindent{\bf Proof.}
In view of \eqref{form(ii)} we have
$$ a_{\ve\beta,\eta}[U]-a_{0,\eta}[U]=2 \ve\mathop{\rm Re}
(\partial_3 U, \dot{\beta} \partial_\vph U)_{\Omega}+\ve^2
\|\dot{\beta}\partial_\vph U\|_\Omega^2\,, $$
and since $|\partial_\vph U|\leq \sup_{x\in\Omega} |x|\: |\nabla
U|$, we get the desired estimate. \hfill$\Box$
\begin{corollary} \label{Cor 6.2} There is an $\ve(\beta)\in (0,1)$ such that for
all $\ve\in(0,\ve(\beta))$ and $\eta\in [-\pi,\pi]$ the estimate
$$
a_{\ve\beta,\eta}[U]\geq \frac12 \|\nabla U\|_\Omega^2
$$
holds. \end{corollary}
\noindent{\bf Proof.}
Observing that $\|\nabla U\|_\Omega^2=a_{0,\eta}[U]$ and using
\eqref{L1} we obtain
$$
\big|a_{\ve\beta,\eta}[U] - \|\nabla U\|_\Omega^2\big|\leq C_\beta\, \ve
\|\nabla U\|_\Omega^2\,,
$$
thus for all $\ve\in (0, C_\beta^{-1})$ the inequality
$$
a_{\ve\beta,\eta}[U]\geq (1-C_\beta\ve)^{-1} \|\nabla U\|_\Omega^2
$$
is valid and it is enough to take $\ve(\beta)= (2C_\beta)^{-1}$.
\hfill$\Box$

\noindent This allows us to estimate the twist effect on the fiber
operator eigenvalues.

\begin{lemma} \label{Lm3.3} To any $k$ there exists a constant $C_{k,\beta}$
such that
$$
|\Lambda_{k,\ve\beta}(\eta)-\Lambda_{k,0}(\eta)|\leq C_{k,\beta}\ve\,.
$$
holds for all $\ve\in (0,\ve(\beta))$ and $\eta\in [-\pi,\pi]$. \end{lemma}
\noindent{\bf Proof.}
Due to the min-max principle, cf.~\cite{BiSo} or \cite{RS4}, we have
$$
\Lambda_{k,\ve\beta}(\eta)=\sup_E\inf_{V\in E\setminus\{0\}}
\frac{a_{\ve\beta,\eta}[V]}{\|V\|_\Omega^2}\,,
$$
where $E$ stands for any subspace in $\cH_\eta$ of codimension
$k-1$. Since the sequence $\{U_{j,0}\}_{j=1}^k$ is chosen
orthonormal in $L^2(\Omega)$ and each $E$ is infinite-dimensional we
can within it an element $U$ of the form
$$
U=\sum_{j=1}^k \alpha_j U_{j,0}\,,\quad \sum_{j=1}^k
|\alpha_j|^2=1\,.
$$
Consequently,
$$
\inf_{V\in E\setminus\{0\}}
\frac{a_{\ve\beta,\eta}[V]}{\|V\|_\Omega^2}\leq
a_{\ve\beta,\eta}[U]\leq a_{0,\eta}[U]+
(a_{\ve\beta,\eta}[U]-a_{0,\eta}[U])\,.
$$
By definition of the form $a_{0,\eta}$ we have
$$
\|\nabla U\|_\Omega^2= a_{0,\eta}[U]=\sum_{j=1}^k
|\alpha_{j}|^2\Lambda_{j,0}(\eta)\leq \Lambda_{k,0}(\eta)\,,
$$
and in combination with inequality \eqref{L1} this yields the
estimate
\begin{equation} \label{L31} \Lambda_{k,\ve\beta}(\eta)\leq
\Lambda_{k,0}(\eta)(1+C_\beta \ve)\,.
\end{equation}
In a similar way one can write
$$
\Lambda_{k,0}(\eta)=\sup_E\inf_{V\in E\setminus\{0\}}
\frac{a_{0,\eta}[V]}{\|V\|_\Omega^2}\,,
$$
where $E$ runs though the same family of subspaces as above.
Repeating the argument we find $U\in E$ of the form $U=\sum_{j=1}^k
\alpha_{j}U_{j,\ve\beta}$ with $\sum_{j=1}^k |\alpha_j|^2=1$ and
infer that
$$
\inf_{V\in E\setminus\{0\}} \frac{a_{0,\eta}[V]}{\|V\|_\Omega^2}\leq
\Lambda_{k,\ve\beta}(\eta)(1+2C_\beta \, \ve)
$$
and
\begin{equation} \label{L32} \Lambda_{k,0}(\eta)\leq
\Lambda_{k,\ve\beta}(\eta)(1+2C_\beta \, \ve)\,.
\end{equation}
The inequalities \eqref{L31} and \eqref{L32} imply the claim of the
Lemma. \hfill$\Box$

\subsection{Asymptotic procedure}

Now we are going to present, in \eqref{4.1} and \eqref{4.2} below,
asymptotic expansions for the eigenvalues
$\Lambda_{k,\ve\beta}(\eta)$ and for the eigenfunctions
$U_{k,\ve\beta}$ as $\ve\to 0$. We consider only simple eigenvalues
$\Lambda_{k,0}(\eta)$ which is sufficient for our purpose; recall
that we deal only with the lowest one in the proof of Theorem~\ref{Th1}. Regarding
the eigenpair $(U_{k,\ve\beta}, \Lambda_{k,\ve\beta}(\eta))$ as a
perturbation of $(U_{k,0}, \Lambda_{k,0}(\eta))$ it is natural to
seek the asymptotic formul{\ae} in the form
\begin{eqnarray} \label{4.1}
&& \Lambda_{k,\ve\beta}(\eta)= \Lambda_{k,0}(\eta)
+\ve \Lambda_k'(\eta)+\ve^2
\Lambda_{k}''(\eta)+\widetilde{\Lambda}_k^\ve(\eta),\\[.3em]
\label{4.2}
&& U_{k,\ve\beta}=U_{k,0}+\ve U_k'+\ve^2 U_k'' +\widetilde{U}_k^\ve,
\end{eqnarray}
where $\Lambda_k'(\eta)$, $\Lambda_{k}''(\eta)$, $ U_k'$, $U_k''$
are the coefficient to be determined and
$\widetilde{\Lambda}_k^\ve(\eta)$, $\widetilde{U}_k^\ve$ the
remainders to be estimated. To begin with, we write the operator
$T_{\ve\beta}(\eta)$ as
\begin{equation} \label{4.3} T_{\ve\beta}(\eta)=T_0(\eta)+\ve T_{1}(\eta)+\ve^2
T_2(\eta), \end{equation}
where $T_0(\eta)$ is the operator associated with the quadratic form
$a_{0,\eta}$, while $T_1(\eta)$ and $T_2(\eta)$ correspond to the
forms
$$
b_{1,\eta}[U]:= 2\mathop{\rm Re} (\partial_3 U, \dot{\beta} \partial_\vph
U)_{\Omega} \quad \mbox{and} \quad b_{2,\eta}[U]:=\|\dot{\beta}
\partial_\vph U\|_\Omega^2
$$
in $\cH_\eta$, respectively. Substituting from
\eqref{4.1}--\eqref{4.3} into the eigenvalue equation
$T_{\ve\beta(\eta)}U_{k,\ve\beta}=\Lambda_{k,\ve\beta}U_{k,\ve\beta}$
and collecting terms of order $\ve$ and $\ve^2$ we get
\begin{eqnarray} \label{4.4} \hspace{-2em}T_0(\eta) U_k'-\Lambda_{k,0}(\eta)
U_k'=\Lambda_k'(\eta) U_{k,0}-T_1(\eta)U_{k,0}\,,
\\[.3em]
\label{4.5} \hspace{-2em}T_0(\eta) U_k''-\Lambda_{k,0}(\eta)
U_k''=\Lambda_k''(\eta)
U_{k,0}-T_2(\eta)U_{k,0}+\Lambda_k'(\eta)U_{k}'-T_1(\eta)U_{k}'\,.
\end{eqnarray}
Since the eigenvalue $\Lambda_{k,0}$ is supposed to be simple and
the resolvent of $T_0(\eta)$ is compact, equation \eqref{4.4} has by
Fredholm alternative one solvability condition, namely its
right-hand side must be orthogonal to $U_{k,0}$ in $L^2(\Omega)$.
Due to the normalization condition we get
\begin{equation} \label{4.6}
\Lambda_k'(\eta)=(T_1(\eta)U_{k,0},U_{k,0})_{\Omega}=2\mathop{\rm Re}
(\partial_3 U_{k,0},\dot{\beta}\partial_\vph U_{k,0})_{\Omega}=0\,,
\end{equation}
where we have used the representation \eqref{10} of the function
$U_{k,0}$ and the following simple observation,
$$
\mathop{\rm Re} (\partial_3 (e^{i\eta x_3} e^{2\pi i m x_3}),
\dot{\beta} e^{i\eta x_3} e^{2\pi i m x_3})_{(-\pi,\pi)}=0\,.
$$
In view of \eqref{4.6} the equations \eqref{4.4} and \eqref{4.5}
simplify to
\begin{eqnarray} \label{4.4new}
&&T_0(\eta) U_k'-\Lambda_{k,0}(\eta)
U_k'=-T_1(\eta)U_{k,0}\,,\\[.3em]
\label{4.5new}
&&T_0(\eta) U_k''-\Lambda_{k,0} (\eta) U_k''=\Lambda_k''(\eta) U_{k,0}-T_2(\eta)U_{k,0}-T_1(\eta)U_{k}'\,.
\end{eqnarray}
Note that the solution $U_k'$ to \eqref{4.4new} is determined up to
a multiple of the eigenfunction $U_{k,0}$, hence without loss of
generality we may assume that
\begin{equation} \label{ort1} (U_k', U_{k,0})_\Omega=0\,. \end{equation}
In the same way, using  the solvability condition of \eqref{4.5new},
we obtain
\begin{eqnarray}
&&\Lambda_k''(\eta)=(T_2(\eta)U_{k,0},U_{k,0})_{\Omega} +2\mathop{\rm Re}
(\partial_3 U_k', \dot{\beta} U_{k,0})_{\Omega}
\nonumber\\[.3em]
&& \quad =\|\dot{\beta}\partial_\vph
\psi_n\|_\Omega^2+\Lambda_{k,0}(\eta)\|U_k'\|_\Omega^2-\|\nabla
U_k'\|_\Omega^2\,.\nonumber \end{eqnarray}
Here we assume that the function $U_{k,0}$ is represented as
$$U_{k,0}(x,\eta)=e^{i\eta x_3}e^{2\pi i m x_3}\psi_n(x')$$
for some $m$ and $n$ and as before we choose $U_k''$ in such a way
that
$$ 
(U_k'',U_{k,0})_\Omega=0\,.
$$ 
Summarizing the above reasoning, we have reached the following
conclusions:
\begin{lemma}
Let $I\subset (-\pi,\pi)$ be a compact set such that
$\Lambda_{k,0}(\eta)$ is simple for all $\eta\in I$. Then for all
$\eta\in I$ the equations \eqref{4.4new} and \eqref{4.5new} have
unique solutions $U_k'$ and $U_k''$ orthogonal to $U_{k,0}$.
Moreover, functions $U_k'$ and $U_k''$ satisfy the estimate
\begin{equation} \label{lm41} \max\{\|U_k'\|_{\cH_\eta}, \|U_k''\|_{\cH_\eta}\}
\leq c_{k,I} \end{equation}
with a constant $c_{k,I}$ independent of $\eta$.
\end{lemma}

Now we are in position to formulate the main result of this section:
\begin{theorem} \label{Th2} Let the eigenvalue $\Lambda_{k,0}(\eta)$ of
operator $T_{0}(\eta)$ be simple for~$\eta$ from a compact set
$I\subset [-\pi,\pi]$, then there exists a $c_{k,\beta,I}$ such that
for all $\ve\in (0,\ve(\beta))$ and $\eta\in I$ one has
$$
\left|\Lambda_{k,\ve\beta}(\eta)-\Lambda_{k,0}(\eta)-\ve^2
\Lambda_k''(\eta)\right|\leq c_{k,\beta,I}\, \ve^3\,.
$$
\end{theorem}

\medskip

Let us first mention that the shortest way to justify the claim of the theorem would be certainly to refer to Sec.~VII.4 of \cite{Ka}. We prefer nevertheless a more explicit discussion because it provides in our view an insight on how the characteristics of the tube vary under the perturbation. In that case a bit of
preliminary work is needed before proceeding to the proof. Let us first introduce a new scalar
product
\begin{equation} \label{scpr} \langle U,V \rangle_{\ve,\eta} :=
A_{\ve\beta,\eta}[U,V] \end{equation}
in the space $\cH_\eta$. According to Lemma~\ref{Lemma 6.1} and
Corollary~\ref{Cor 6.2} the corresponding norm is uniformly
equivalent to the gradient norm for all $\ve\in (0,\ve(\beta))$. The
Hilbert space with the inner product \eqref{scpr} will be denoted by
$\cH_{\ve,\eta}$. Next we define the operator $L_{\ve,\eta}$ by the
formula
$$
\langle L_{\ve,\eta} U, V\rangle_{\ve,\eta}=(U,V)_\Omega\,;
$$
it is easy to check that it is compact and self-adjoint and its
discrete spectrum consists of the eigenvalue sequence
\begin{equation} \label{nu} \nu_{k,\ve}(\eta)=
(\Lambda_{k,\ve\beta}(\eta))^{-1}\,. \end{equation}

\noindent \emph{Proof of Theorem~\ref{Th2}:} According to
Lemma~\ref{Lm3.3} it is sufficient to prove that there is a constant
$c_{k,\beta,I}>0$ such that the interval
$$(\Lambda_{k,0}(\eta)+\ve^2 \Lambda_k''(\eta)-c_{k,\beta,I} \ve^3,
\Lambda_{k,0}(\eta)+\ve^2 \Lambda_k''(\eta)+c_{k,\beta,I} \ve^3)$$
contains at least one member of the sequence
$\{\Lambda_{j,\ve\beta}(\eta)\}_{j\geq 1}$. Equivalently, it is
enough to establish the existence of at least one eigenvalue from
the sequence \eqref{nu} in the interval
$$
((\Lambda_{k,0}(\eta)+\ve^2 \Lambda_k''(\eta))^{-1}-c'_{k,\beta,I} \ve^3,
(\Lambda_{k,0}(\eta)+\ve^2 \Lambda_k''(\eta))^{-1}+c'_{k,\beta,I} \ve^3)
$$
for some $c'_{k,\beta,I}>0$. To do this we construct a function
$W\in \cH_{\ve,\eta}\setminus \{0\}$ such that
$$
\|L_{\ve,\eta} W-(\Lambda_{k,0}(\eta)+\ve^2
\Lambda_k''(\eta))^{-1}W\|_{\cH_{\ve,\eta}} \leq c'_{k,\beta,I}
\ve^3 \|W\|_{\cH_{\ve,\eta}}\,,
$$
namely, we put $W:=U_{k,0}+\ve U_k'+\ve^2 U_k''$. Then
\begin{equation} \|(L_{\ve,\eta}-(\Lambda_{k,0}(\eta)+\ve^2
\Lambda_k''(\eta))^{-1})W\|_{\cH_{\ve,\eta}} = \sup_{V}\langle(
L_{\ve,\eta} -(\Lambda_{k,0}(\eta)+\ve^2 \Lambda_k''(\eta))^{-1})W,
V\rangle_{\ve,\eta}\,, \nonumber \end{equation}
where $V\in \cH_{\ve,\eta}$ runs over over unit-length vectors,
\begin{equation} \label{V} \|V\|_{\cH_{\ve,\eta}}=1\,. \end{equation}
Let us observe the following chain of equalities
\begin{align*}
\tau:=&\langle( L_{\ve,\eta} -(\Lambda_{k,0}(\eta)+\ve^2
\Lambda_k''(\eta))^{-1})W,
V\rangle_{\ve,\eta}\\[.3em]
=&(W,V)_\Omega - (\Lambda_{k,0}(\eta)+\ve^2 \Lambda_k''(\eta))^{-1}
A_{\ve\beta,\eta}[W,V]=\\[.3em]
=&(W,V)_\Omega - (\Lambda_{k,0}(\eta)+\ve^2
\Lambda_k''(\eta))^{-1}((T_0(\eta)+\ve T_1(\eta)+\ve^2
T_2(\eta))W,V)_\Omega\\[.3em]
=&(\Lambda_{k,0}(\eta)+\ve^2 \Lambda_k''(\eta))^{-1} \\[.3em] &\;\times (\ve^3
\Lambda''(\eta)U'_k+\ve^4\Lambda''(\eta)U_k''-\ve^3T_2U'-\ve^3T_1U''-\ve^4
T_2U'',V)_\Omega\,.
\end{align*}
Since the expression $(\Lambda_{k,0}(\eta)+\ve^2
\Lambda_k''(\eta))^{-1}$ is bounded from above on $I$, in view of
the estimates \eqref{lm41} and normalization condition \eqref{V} we
can conclude that there is a $C_{k,\beta,I}$ such that
$$
|\tau|\leq C_{k,\beta,I}\,\ve^3
$$
holds for all $\ve$ small enough. To complete the proof it is
sufficient to observe that by virtue of the estimates \eqref{lm41}
combined with Lemma~\ref{Lemma 6.1} and Corollary~\ref{Cor 6.2} the
expression $\|W\|_{\cH_{\ve,\eta}}$ is uniformly bounded away from
zero for $\eta\in I$ and $\ve\in (0,\ve(\beta))$.
\hfill{$\Box$}
\begin{corollary} \label{Cor2} There is a $\bfc_\beta>0$ such that for all
$\eta\in [-\pi/2,\pi/2]$ we have
$$
|\Lambda_{1,\ve\beta}(\eta)-\lambda_1-\eta^2|\leq \bfc_\beta \ve^2\,.
$$
\end{corollary}

\subsection{Concluding the proof}

We now able to prove the main result of this section,
Theorem~\ref{Th1}. Let us divide the interval $[-\pi,\pi]$ into
three parts. On the first part, $I_1:=\{\eta \colon |\eta|>\pi/2\}$,
we write using Lemma~\ref{Lm3.3}
$$
\Lambda_{1,\ve\beta}(\eta)\geq \Lambda_{1,0}(\eta)-C_{1,\ve}\ve
=\lambda_1+\eta^2-C_{1,\ve}\ve\geq \lambda_1+(\pi/2)^2-C_{1,\ve}\ve\,.
$$
On the second part, $I_2:=\{\eta \colon |\eta|\in [\alpha \ve,
\pi/2]\}$, we use Corollary~\ref{Cor2} to obtain
$$
\Lambda_{1,\ve\beta}(\eta)\geq \lambda_1+\eta^2-\bfc_\beta \ve^2\geq
\lambda_1 +(\alpha^2-\bfc_\beta)\ve^2\,;
$$
the parameter $\alpha$ will be specified below. As for the third
part, $I_3:=\{\eta \colon |\eta|<\alpha \ve\}$, applying the
estimate from Theorem~\ref{Th2} we get
$$
\Lambda_{1,\ve\beta}(\eta)\geq \lambda_1+\eta^2+\ve^2 \Lambda''_1(\eta)-c_{1,\beta,I_3} \ve^3\geq \lambda_1 +\ve^2\Lambda''_1(\eta)-c_{1,\beta,I_3} \ve^3\,,
$$
where $\Lambda''_1(\eta)$ can be calculated by the formula
\begin{align*}
\Lambda_1''(\eta)&=
\|\dot{\beta}\partial_\vph \psi_1\|_\Omega^2 +\Lambda_{1,0}(\eta)\|U_1'\|_\Omega^2-\|\nabla
U_1'\|_\Omega^2\\[.3em]
&\geq\int_0^1 |\dot{\beta}(x_3)|^2\,\mathrm{d}x_3 \int_\omega
|\partial_\vph \psi_1(x')|^2\,\mathrm{d}x' - C \alpha^2 \ve^2\,.
\end{align*}
The last inequality holds true, because the solution of the equation
\eqref{4.4new} for $k=1$ with orthogonality condition \eqref{ort1}
satisfies the estimate
\begin{align*}
\|U'_1\|_\Omega&+\|\nabla U'_1\|_\Omega\leq
C\|T_1(\eta)U_{1,0}\|_\Omega \\[.3em] &= C
|\eta|\|\dot{\beta}\|_{L^2(0,1)}\|\partial_\vph
\psi_1\|_{L^2(\omega)}\,.
\end{align*}
Now we put $\alpha^2:=\bfc_\beta+\int_0^1
|\dot{\beta}(x_3)|^2\,\mathrm{d}x_3 \int_\omega |\partial_\vph
\psi_1(x')|^2\,\mathrm{d}x\,\mathrm{d}y$ and obtain for all
sufficiently small $\ve$ the estimate
$$
\lambda_\dagger(\ve\beta)\geq \lambda_1 +\ve^2 \int_0^1
|\dot{\beta}(x_3)|^2\,\mathrm{d}x_3 \int_\omega |\partial_\vph
\psi_1(x')|^2\,\mathrm{d}x' -C_\beta \ve^3\,.
$$
Together with the result of Theorem~\ref{Th2} for $\eta=0$ and $k=1$
we get in this way the asymptotic expansion
$$
\lambda_\dagger(\ve\beta)= \lambda_1 +\ve^2 \int_0^1
|\dot{\beta}(x_3)|^2\,\mathrm{d}x_3 \int_\omega |\partial_\vph
\psi_1(x')|^2\,\mathrm{d}x' +O(\ve^3)\,.
$$
The claim of Theorem~\ref{Th1} is then a simple consequence of this
formula and relation~\eqref{assump}. \hfill{$\Box$}

\medskip

As the last remark, one is naturally interested whether the property expressed by Theorem~\ref{Th1} could remain valid beyond the perturbative regime. This question would require a different approach and one cannot exclude that such a claim may not hold generally.

\subsection*{Acknowledgements}

In its first part the work was supported by the project 17-01706S of the Czech Science
Foundation (GA\v{C}R). The results of the second part of the paper, notably Theorems~\ref{Th1} and~\ref{Th2}, 
were obtained with the support of the Russian Science Foundation Grant 14-21-00035.

\end{document}